\documentclass[12pt,twoside]{article}

\usepackage{epsfig,amsmath,amssymb}
\usepackage{amsfonts}
\usepackage{amssymb}
\usepackage{amssymb}
\usepackage{color}
\usepackage{enumerate}

\date{}

\begin{document}


\centerline{}

\centerline{}

\centerline {\Large{\bf Inscribed triangles in the unit sphere}}
\centerline {\Large{\bf  and a new class of geometric constants}}

\centerline{}

\centerline{\bf {}}

\centerline{Bingren Chen, Qi Liu, Yongjin Li$^*$}

\centerline{Department of Mathematics,
      Sun Yat-sen University}

\centerline{ Guangzhou, 510275, P. R. China}

\let\thefootnote\relax
\footnotetext{$^*$Corresponding author. E-mail: stslyj@mail.sysu.edu.cn}

\newtheorem{Theorem}{\qquad Theorem}[section]

\newtheorem{Definition}[Theorem]{\qquad Definition}

\newtheorem{Proposition}[Theorem]{\qquad Proposition}

\newtheorem{Corollary}[Theorem]{\qquad Corollary}

\newtheorem{Proof}[Theorem]{\qquad Proof}

\newtheorem{Lemma}[Theorem]{\qquad Lemma}

\newtheorem{Remark}[Theorem]{\qquad Remark}

\newtheorem{Example}[Theorem]{\qquad Example}

\newtheorem{Question}[Theorem]{\qquad Question}

\begin{abstract}
We will introduce a new geometric constant $G_L(X)$ based on the constant $H(X)$ proposed by Gao. We first further survey the constant $H(X)$ and discuss some of the properties of this constant that have not yet been discovered. Next, we focus on a new constant $G_L(X)$ along with some of its basic properties. In addition, we show some relations between the well-known geometric constants and $G_L(X)$ through some inequalities. Finally, we characterize some generalized forms of the constant $G_L(X)$.
\end{abstract}



\section{Introduction}
\qquad Geometric constant has got widespread attention, because it can not only essentially reflect the geometric properties of a space $X$, but also help us to quantitatively study the space. As a research tool, geometric constants are also of great significance for their own research, such as the estimation of some specific space constants, the relationship between some constants, and the relationship between the constant value of the space itself and the value of the dual space.

 The most classic constants are $C_{\rm{NJ}}(X)$ and $J(X)$. Among them, $J(X)$ was proposed by Gao  et al. in 1991. It was put forward on the basis of James’ characterization of uniformly non-square space.
A further research has been conducted on whether the space has normal structure. For the result and application of this constant, please refer to the article \cite{Mizuguchi2021,Zuo2020,
	Kato2001,Ahmad,LQ}.

The constant $H(X)$ was proposed by Gao  in 2000. Unlike the existing constants, this constant combines the unit sphere inscribed with an equilateral triangle and innovatively infuses the quantity
$ 2x-y$. For its geometric significance, we will give an explanation in the subsequent part of this article. It is worth mentioning that the introduction of the quantity $2x-y$ here has actually increased the difficulty of the study to a large extent. Some classical geometric properties and the connection of constants will be invalid. On the one hand, unlike the constant $J(X)$, the constant $H(X)$ is no longer symmetrical.
Most of the geometric constants we know all have a symmetrical structure. On the other hand, there are also differences in some techniques and inequalities scaling when dealing with asymmetric geometric constants.

\section {Preliminaries}
\qquad We will assume throughout the paper that $X$ represents a non-trivial Banach space, that is, $\operatorname{dim} X \geq 2$, and use $S_X$ and $B_X$ to represent the unit sphere and unit ball of $X$, respectively.

\begin{Definition}\cite{James1964} 
	A Banach space $X$ is called uniformly convex, if for any $\varepsilon>0$, there exists $\delta>0$, such that for any $x, y\in S(X)$ with $\Vert x-y\Vert>\varepsilon, \Vert (x+y)/2\Vert<1-\delta$.
\end{Definition}

The Clarkson modulus of convexity of a Banach space $X$
is defined as follows \cite{JAC}:
$$\delta_X(\epsilon)=\inf \bigg\{1-\frac{\Vert x+y\Vert}{2}:x,y\in S_X, \Vert x-y\Vert\geqslant \epsilon \bigg\}.$$

The constant $J(X)$
$$
J(X)=\sup \left\{\min (\|x+y\|,\|x-y\|): x, y \in S_{X}\right\}
$$
is called the James constant \cite{Gao1990}, and the  constant $C_{\rm{NJ}}(X)$  \cite{Clarkson1937} is introduced by Clarkson to describe the inner product space.
In a sense, it can be understood as the following formula, for which
$$
C_{\mathrm{NJ}}(X)=\sup \left\{\frac{\|x+y\|^{2}+\|x-y\|^{2}}{2\left(\|x\|^{2}+\|y\|^{2}\right)}: (x,y)\neq (0,0)\right\}.
$$

The  constant $C_Z(X)$ was introduced by G. Zbăganu \cite{ZG}:
$$
C_{\mathrm{Z}}(X)=\sup \left\{\frac{\|x+y\|\|x-y\|}{\|x\|^{2}+\|y\|^{2}}:x, y \in X, (x, y) \neq(0,0)\right\}.
$$
$C_Z(X)$ and  $C_{\mathrm{NJ}}(X)$ seem to be compatible.
 It is worth mentioning that  Alonso and Martin \cite{JA3}  gave  a counterexample
that $ C_Z(X)\neq C_{\mathrm{NJ}}(X)$.

The next constant, which is closely related to the new constant we studied, was defined by Gao. The constant $H(X)$ is defined as \cite{Gao2000}
$$H(X)=\sup\{\min \{\Vert x+y\Vert,\Vert 2x-y\Vert\}:x, y, x-y\in S(X)\}.$$

Listed below are some of the results of the constant $H(X)$ given in \cite{Gao2000}:

(i) For a Banach space, $H(X)\leq 2$.

(ii) If $X$ is a Hilbert space, then $H(X)= \sqrt{3}$.

(iii) If $H(X)< 2$, then $X$ is uniformly non-square.

Next, we give more properties of the $H(X)$ constant, because some of the techniques used in it only appeared after the $H(X)$ constant was introduced.
The proof of Theorem \ref{Theorem1.1} is based on the method published in 2001 by Kato et al.

\begin{Theorem}
	Let Banach space $X$ be finite-dimensional, if $H(X)=2$, then $X$ is not strictly convex.
\end{Theorem}
\begin{Proof}
Assume that $H(X)=2$. Since the unit sphere of finite-dimensional Banach space is compact, so there exist $x_0,y_0,x_0-y_0\in S(X)$, such that
$\Vert x_0+y_0\Vert=2$, which implies that point $\frac{x_0+y_0}{2}$ is not extreme point of closed unit sphere of $X$, so $X$ is not strictly convex.
\end{Proof}

\begin{Theorem}\label{Theorem1.1}
	Let $X$ be a Banach space, then
	$$2H(X)-2\leq H(X^{*})\leq\frac{1}{2}H(X)+1.   $$
\end{Theorem}

\begin{Proof}
For any $x,y,x-y\in S(X)$, we have $\Vert x+y\Vert \leq\Vert x\Vert+\Vert y\Vert=2$, $\Vert 2x-y\Vert\leq 2$, thus
$$\Vert x+y\Vert+\Vert 2x-y\Vert\leq\sup\{\min\{\Vert x+y\Vert,\Vert 2x-y\Vert\}\}+2.$$
For any $\varepsilon>0$, there exist $x,y,x-y\in S(X)$, such that
$$\min\{\Vert x+y\Vert,\Vert 2x-y\Vert\}\geq H(X)-\varepsilon.$$
Additionally, by Hahn-Banach Theorem, there exist functionals $u^{*}, v^{*}\in S(X^{*})$, such that
$$u^{*}(x+y)=\Vert x+y\Vert, v^{*}(y-2x)=\Vert y-2x\Vert=\Vert 2x-y\Vert.$$
Thus
$$
\begin{aligned}
	H(X^{*})&= \sup\{\Vert u^{*}+v^{*}\Vert,\Vert 2u^{*}-v^{*}\Vert\}\\
	&\geq\Vert  u^{*}+v^{*}\Vert+\Vert 2u^{*}-v^{*}\Vert-2\\
	&\geq(u^{*}+v^{*})(y-x)+(2u^{*}-v^{*})(x)-2\\
	&=u^{*}(x+y)+v^{*}(y-2x)-2\\
	&=\Vert x+y\Vert+\Vert 2x-y\Vert-2\\
	&\geq 2(\min\{\Vert x+y\Vert,\Vert 2x-y\Vert\})-2\\
	&\geq 2(H(X)-\varepsilon)-2.
\end{aligned}
$$
Since $\varepsilon$ can be arbitrarily small, so
$$H(X^{*})\geq 2H(X)-2.$$
To prove the right side of inequality, assume that $u^{*}, v^{*}\in S(X^{*})$, there exist $x, y, x-y\in S(X)$, such that
$$(u^{*}+v^{*})(y-x)>\Vert u^{*}+v^{*}\Vert-\varepsilon, (2u^{*}-v^{*})(x)>\Vert 2u^{*}-v^{*}\Vert-\varepsilon.$$
Therefore,
$$
\begin{aligned}
	\sup\{\min\{\Vert u^{*}+v^{*}\Vert,\Vert 2u^{*}-v^{*}\Vert\}\}&\leq\frac{1}{2}(\Vert u^{*}+v^{*}\Vert+\Vert 2u^{*}-v^{*}\Vert)\\
	&\leq\frac{1}{2}[(u^{*}+v^{*})(y-x)+(2u^{*}-v^{*})(x)+2\varepsilon]\\
	&=\frac{1}{2}[u^{*}(x+y)+v^{*}(y-2x)+2\varepsilon]\\
	&\leq\frac{1}{2}(\Vert u^{*}\Vert\cdot\Vert x+y\Vert+\Vert v^{*}\Vert\cdot\Vert 2x-y\Vert)+2\varepsilon\\
	&=\frac{1}{2}(\Vert x+y\Vert+\Vert 2x-y\Vert)+2\varepsilon\\
	&\leq\frac{1}{2}\sup\{\min\{\Vert x+y\Vert,\Vert 2x-y\Vert\}\}+2+2\varepsilon\\
	&\leq\frac{1}{2} H(X)+2+2\varepsilon.
\end{aligned}
$$
Since $\varepsilon$ can be arbitrarily small, so we prove the result as desired.
\end{Proof}

\section{The constant $G_L(X)$}
\qquad We introduce a new constant based on the constant $H(X)$. We begin by introducing the following key definition:

$$G_L(X)=\sup\{\Vert x+y\Vert^{2}+\Vert 2x-y\Vert^{2}~:~ \|x\|=\|y\|=\|x-y\|=1\}.$$

\begin{figure}[ht]
	\centering
	\includegraphics[scale=0.6]{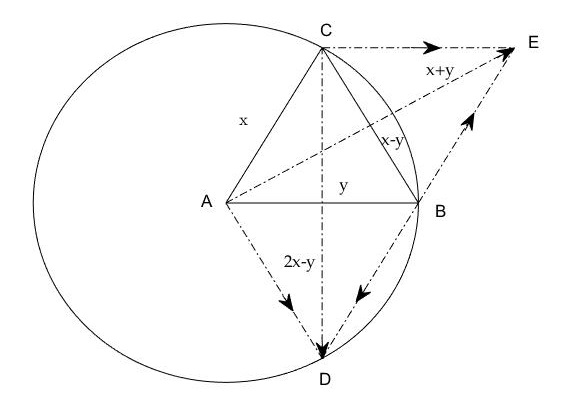}
	\caption{ Geometric explanation on the unit sphere}\label{fig:simu}
	\label{fig:1}{}
\end{figure}

The geometric background of $G_{L}(X)$ is shown in Figure 1: consider the unit sphere on the Euclidean plane with $\Vert x\Vert=\Vert y\Vert=\Vert y-x\Vert=1$. Assume that $\mathop{AC}\limits^{\rightarrow}=x$, $\mathop{AB}\limits^{\rightarrow}=y$, then $\mathop{AE}\limits^{\rightarrow}=x+y$. Assume that $\mathop{CB}\limits^{\rightarrow}=y-x$, then $\mathop{DC}\limits^{\rightarrow}=2x-y$. Apparently, triangle $\triangle ABC$ is equilateral. Moreover, $\mathop{AD}\limits^{\rightarrow}=x-y$, $\mathop{CE}\limits^{\rightarrow}=y$, $\mathop{DE}\limits^{\rightarrow}=2x$.

Since the Clarkson type inequalities
$$2(\Vert x\Vert ^{p}_{p}+\Vert y\Vert ^{p}_{p})\leq\Vert x+y\Vert ^{p}_{p}+\Vert x-y\Vert ^{p}_{p}\leq 2^{p-1}(\Vert x\Vert ^{p}_{p}+\Vert y\Vert ^{p}_{p})$$
holds for $p\geq 2$, $x,y$ in $\ell_p$ , we take $p=2$, where the demonstration of the given Figure 1 hinges on, then we have
$$\begin{aligned}&2\bigg(\bigg\Vert \frac{(2x-y)+y}{2}\bigg\Vert^{2}+\bigg\Vert\frac{(2x-y)-y}{2}\bigg\Vert^{2}\bigg)\\&\leq\Vert 2x-y\Vert^{2}+\Vert y\Vert^{2}\\&\leq2^{2-1}\bigg(\bigg\Vert \frac{(2x-y)+y}{2}\bigg\Vert^{2}+\bigg\Vert\frac{(2x-y)-y}{2}\bigg\Vert^{2}\bigg),
\end{aligned}$$
i.e.
$$\begin{aligned}2(\Vert x\Vert^{2}+\Vert x-y\Vert^{2})&\leq\Vert 2x-y\Vert^{2}+\Vert y\Vert^{2}\\&\leq 2(\Vert x\Vert^{2}+\Vert x-y\Vert^{2}),
	\end{aligned}$$
thus
$$2(\Vert x\Vert^{2}+\Vert x-y\Vert^{2})=\Vert 2x-y\Vert^{2}+\Vert y\Vert^{2},$$
therefore,
$${\mathop{CD}\limits^{\rightarrow}}^{2}+{\mathop{AB}\limits^{\rightarrow}}^{2}=2({\mathop{BD}\limits^{\rightarrow}}^{2}+{\mathop{AD}\limits^{\rightarrow}}^{2})$$
holds, and then the quadrilateral $\Box ACBD$ satisfies the parallelogram law.
Similarly, we can deduce that the quadrilateral $\Box ABEC$ satisfies the parallelogram law, either.
In general, the quadrilateral $\Box ACED$ can be described as half a Hexagon.

\begin{Theorem}
	For any Banach space $X$, we have $\frac{9}{2}\leq G_L(X) \leq 8$.
\end{Theorem}

\begin{Proof}
Since
$$\begin{aligned}
	\Vert x+y\Vert+\Vert 2x-y\Vert&\geq \Vert x+y+2x-y\Vert
	\\&=\Vert 3x\Vert\\&=3
\end{aligned}$$
and hence
$$\Vert 2x-y\Vert\geq 3-\Vert x+y\Vert.$$
Then we can deduce that
$$\begin{aligned}
	\Vert x+y\Vert^2+\Vert 2x-y\Vert^2&\geq \Vert x+y\Vert^2+(3-\Vert x+y\Vert)^2
	\\&=2\Vert x+y\Vert^2+9-6\Vert x+y\Vert\\&=
	2\bigg[(\Vert x+y\Vert-\frac{3}{2})^2+\frac{9}{4}\bigg]
	\\&\geq \frac{9}{2}.
\end{aligned}$$
\end{Proof}

On the other hand, by $\|2x-y\|\leq \|x\|+\|x-y\|=2$,
which implies that $G_L(X)\leq 8$, as desired.

\begin{Corollary}
Let Banach space $X$ be finite dimensional, if $G_L(X)=8$, then $X$ is not strictly convex.
\end{Corollary}

\begin{Example}
	Let $X$ be $\ell_{\infty}$ endowed with the norm $\Vert x\Vert_{\infty}=\sup_{n}\vert x_{n}\vert$ for $(x_{n})\in\ell_{\infty}$. Then $G_L(\ell_{\infty})=8$.
\end{Example}
Assume that $x_0=(1,1,0,\ldots), y_0=(1,0,0,\dots)\in S(\ell_{\infty})$, then $x_0-y_0=(0,1,0,\ldots)\in S(\ell_{\infty})$.
We can get $$\begin{aligned}G_L(\ell_{\infty})&\geq\Vert x_0+y_0\Vert^{2}+\Vert 2x_0-y_0\Vert^{2}\\&=2^{2}+2^{2}=8,
	\end{aligned}$$
 which implies that $G_L(\ell_{\infty})= 8$.

\begin{Example}
	Let $X$ be $\ell_p, 1<p<\infty$, then
	$G_L(X)\leq 2(2^p-1)^{\frac{2}{p}}$ for $p\geq 2$ and
	$G_L(X)\leq 2\cdot3^{\frac{2}{p}}$ for $1\leq p< 2$.
\end{Example}
Applying Clarkson inequality that when
$p\geq 2, x,y\in X$, we have
$$\begin{aligned}2(\Vert x\Vert^p+\Vert y\Vert^p)
	&\leq \Vert x+y\Vert^p+\Vert x-y\Vert^p
	\\&\leq (\Vert x\Vert+\Vert y\Vert)^p+
	|\Vert x\Vert-\Vert y\Vert|^p,
\end{aligned}$$
when
$1\leq p< 2, x,y\in X$, we have
$$\begin{aligned}2(\Vert x\Vert^p+\Vert y\Vert^p)
	&\geq \Vert x+y\Vert^p+\Vert x-y\Vert^p
	\\&\geq (\Vert x\Vert+\Vert y\Vert)^p+
	|\Vert x\Vert-\Vert y\Vert|^p.
\end{aligned}$$

Then we can deduce that
when
$p\geq 2, x,y,x-y\in S_X$,
$$\begin{aligned}\Vert x+y\Vert^p&\leq (\Vert x\Vert+\Vert y\Vert)^p+
	|\Vert x\Vert-\Vert y\Vert|^p-\Vert x-y\Vert^p\\
	&= 2^p-1,\end{aligned}$$
$$\begin{aligned}
	\Vert 2x-y\Vert^p&\leq(\Vert x\Vert+\Vert x-y\Vert)^p+
	|\Vert x\Vert-\Vert x-y\Vert|^p-\Vert x-(x-y)\Vert^p\\&=  2^p-1,\end{aligned}$$
when
$1\leq p< 2, x,y,x-y\in S_X$,
$$\begin{aligned}\Vert x+y\Vert^p&\leq2(\Vert x\Vert^p+\Vert y\Vert^p)-\Vert x-y\Vert^p\\&\leq 2(\Vert x\Vert^2+\Vert y\Vert^2)-\Vert x-y\Vert^2\\&=3,\end{aligned}$$
$$\begin{aligned}\Vert 2x-y\Vert^p&\leq2(\Vert x\Vert^p+\Vert x-y\Vert^p)-\Vert x-(x-y)\Vert^p\\&\leq 2(\Vert x\Vert^2+\Vert x-y\Vert^2)-\Vert x-(x-y)\Vert^2\\&= 3,\end{aligned}$$
as desired.

\begin{Proposition}
	If $X$ is an inner product space, then $G_L(X)=6$.
\end{Proposition}
\begin{Proof}
For any $x, y, x-y\in S(X)$, we have
$$\Vert x+y\Vert^{2}=2(\Vert x\Vert^{2}+\Vert y\Vert^{2})-\Vert x-y\Vert^{2}=3,$$
$$\Vert 2x-y\Vert^{2}=2(\Vert x\Vert^{2}+\Vert x-y\Vert^{2})-\Vert x-(x-y)\Vert^{2}=3,$$
which implies that  $G_L(X)=6$.
\end{Proof}


\begin{Proposition}
	For any Banach space $X$, we have
	$$4(1-\delta_{X}(1))^2\geq G_L(X)-4.$$
\end{Proposition}
\begin{Proof}
 First, note that
$$
\begin{aligned}\delta_{X}(t)&=\inf \left\{1-\frac{1}{2}\|x+y\|:\|x\|=\|y\|=1,\|x-y\| \geqslant \varepsilon\right\} \\& =\inf \left\{1-\frac{1}{2}\|x+y\|:\|x\|=\|y\|=1,\|x-y\|=\varepsilon\right\}.\end{aligned}
$$
Then we can deduce that
$$
\delta_{X}(\|x-y\|) \leqslant 1-\frac{\|x+y\|}{2}
$$
for any $\|x\|=\|y\|=\|x-y\|=1$.

Applying the triangle inequality, we have the following inequality estimate:
$$\begin{aligned}
	(1-\delta_{X}(\|x-y\|))^2+4&\geq \frac{\|x+y\|^2}{4}+1+3
	\\&=\frac{1}{4}(\Vert x+y\Vert^{2}+\|x\|^2 +\Vert x-y\Vert^{2}+2\|x\|\Vert x-y\Vert)+3
	\\&\geq \frac{1}{4}(\Vert x+y\Vert^{2}+\Vert 2x-y\Vert^{2})+3
\end{aligned}$$
for any $\|x\|=\|y\|=\|x-y\|=1$,
and hence
$$(1-\delta_{X}(\|x-y\|))^2+4\geq \frac{1}{4}(\Vert x+y\Vert^{2}+\Vert 2x-y\Vert^{2})+3$$
for any $\|x\|=\|y\|=\|x-y\|=1$.
This means that
$$4(1-\delta_{X}(1))^2\geq G_L(X)-4,$$
as desired.
\end{Proof}

\begin{Proposition}\label{Proposition1}
	For any Banach space $X$, we have
	$$H(X)^2\leq \frac{1}{2}G_L(X).$$
\end{Proposition}
\begin{Proof} Consider  points $x$ and  $y$ in space $X$, we  have the following estimates
$$\begin{aligned}\min\{\Vert x+y\Vert^2,\Vert 2x-y\Vert^2\}&\leq \frac{(\Vert x+y\Vert+\Vert 2x-y\Vert)^2}{4}\\&\leq \frac{\Vert x+y\Vert^2+\Vert 2x-y\Vert^2}{2},
\end{aligned}$$
as desired.
\end{Proof}

\begin{Theorem}
	Let $X$ be a Banach space, we have
	(i)$\Rightarrow$(ii) $\Rightarrow$(iii).

	(i) $G_L(X)<8$.

	(ii) $H(X)<2$.

	(iii) $X$ is uniformly non-square.
\end{Theorem}
\begin{Proof}
(i)$\Rightarrow$ (ii). Applying Proposition  \ref{Proposition1}.

(ii)$\Rightarrow$ (iii).  Using the result in [\cite{Gao2000}, p. 243, Theorem 2.10], if $H(X)< 2$, then $X$ is uniformly non-square.
\end{Proof}

\begin{Corollary}
	If $X$ is not super-reflexive, then $G_L(X)=8$.
\end{Corollary}

\section{The constant $G_L(X,p)$ and $C_L(X)$}

\qquad Next, we will  introduce two generalizations of constants and consider their related properties.

Given any Banach space $X$ and a number $p\in[1,\infty)$, another geometric constant $G_L(X,p)$ is defined by
$$G_L(X,p)=\sup\{\Vert x+y\Vert^{p}+\Vert 2x-y\Vert^{p}~:~ \|x\|=\|y\|=\|x-y\|=1\}.$$

\begin{Proposition}
	For any Banach space $X$, we have $2^{1-p}\cdot3^{p}\leq G_L(X,p) \leq2^{p+1}$.
\end{Proposition}
\begin{Proof}
By the convexity of the function $f(u)=u^{p}$ on $[0,\infty)$, we have the following inequality:

$$\frac{\Vert x+y\Vert^{p}+\Vert 2x-y\Vert^{p}}{2}\geq\bigg(\frac{\Vert x+y\Vert+\Vert 2x-y\Vert}{2}\bigg)^{p}.$$
Then we can deduce that
$$\begin{aligned}
\frac{\Vert x+y\Vert^{p}+\Vert 2x-y\Vert^{p}}{2}&\geq(\frac{\Vert x+y\Vert+\Vert 2x-y\Vert}{2})^{p}\\&\geq(\frac{\Vert x+y+2x-y\Vert}{2})^{p}\\&=2^{-p}\cdot3^{p},
\end{aligned}
$$
hence
$$\Vert x+y\Vert^{p}+\Vert 2x-y\Vert^{p}\geq 2^{1-p}\cdot3^{p}.$$

On the other hand,  according to the triangle inequality, we can have the following estimates
$$\begin{aligned}
	\Vert x+y\Vert^{p}+\Vert 2x-y\Vert^{p}&\leq (\Vert x\Vert+\Vert y\Vert)^{p}+(\Vert x\Vert+\Vert x-y\Vert)^{p}
	\\&=2^{p+1}.
\end{aligned}$$
This completes the proof.
\end{Proof}

\begin{Theorem}
	For any $1<p<\infty$ and any Banach space $X$, the following inequality holds:
	$$H(X)\leq2^{-\frac{1}{p}}\sqrt[p]{G_{L}(X,p)}.$$
\end{Theorem}
\begin{Proof}
Indeed, if $1<p<\infty$, then for any $x,y,x-y\in S(X)$, we have
$$2(\min\{\Vert x+y\Vert, \Vert 2x-y\Vert\})^{p}\leq \Vert x+y\Vert^{p}+\Vert 2x-y\Vert^{p},$$
which implies that
 $$\min\{\Vert x+y\Vert, \Vert 2x-y\Vert\}\leq2^{-\frac{1}{p}}\sqrt[p]{G_{L}(X,p)},$$
and the proof is completed.

\end{Proof}

Combining our new constant $G_L(X)$ with the constant  $C_{\mathrm{Z}}(X)$, we define the $C_{L}(X)$ constant as follow:
$$C_{L}(X)=\sup\{\Vert x+y\Vert\cdot\Vert 2x-y\Vert: \Vert x\Vert=\Vert y\Vert =\Vert x-y\Vert=1\}.$$

\begin{Proposition}\label{Proposition3}
	Let $X$ be a Banach space, then
	$$C_{L}(X)\leq\frac{1}{2}G_{L}(X).$$
\end{Proposition}
\begin{Proof}
Consider  points $x$ and  $y$ in space $X$, we  have the following estimates
$$2\Vert x+y\Vert\cdot\Vert 2x-y\Vert\leq\Vert x+y\Vert^{2}+\Vert 2x-y\Vert^{2},$$
as desired.
\end{Proof}

\begin{Example}
	Let $X$ be $\ell_{\infty}$ endowed with the norm  $\Vert x\Vert_{\infty}=\sup_{n}\vert x_{n}\vert$  for $(x_{n})\in\ell_{\infty}$. Then $C_L(\ell_{\infty})=4$.
\end{Example}
Assume that $x=(1,1,0,\ldots), y=(1,0,0,\dots)\in S(\ell_{\infty})$, then $x-y=(0,1,0,\ldots)\in S(\ell_{\infty})$.
Then $C_L(\ell_{\infty})\geq\Vert x+y\Vert\cdot\Vert 2x-y\Vert=2\cdot2=4$.
Since  $C_L(\ell_{\infty})\leq \frac{1}{2}G_{L}(\ell_{\infty})=4$, so $C_L(\ell_{\infty})= 4$.

\begin{Proposition}
	For any Banach space $X$, we have
	$$C_{L}(X)\geq \frac{1}{2}(9-G_L(X)).$$
\end{Proposition}
\begin{Proof}
Since
$$\begin{aligned}
	\Vert x+y\Vert+\Vert 2x-y\Vert&\geq \Vert x+y+2x-y\Vert
	\\&=\Vert 3x\Vert\\&=3
\end{aligned}$$
and hence
$$\Vert 2x-y\Vert^2+\Vert x+y\Vert^2+2\Vert 2x-y\Vert\Vert x+y\Vert\geq 9.$$
Then we can deduce that
$$
\Vert x+y\Vert\Vert 2x-y\Vert\geq \frac{1}{2}(9-G_L(X)),$$
which implies that $C_{L}(X)\geq \frac{1}{2}(9-G_L(X))$, as desired.
\end{Proof}

\section{Conclusions}
\qquad {We introduce a new geometric constant $G_L(X)$, which is based on the equilateral triangle in the unit sphere and is closely related to the constant which is introduced by Gao. We give the relation between the famous geometric constant and $G_L(X)$ through some inequalities and use it to describe some geometric properties. At the same time, we also study the generalized form of $G_L(X)$. The geometric background of the constant defined is very intuitive. There are many similarities with the classical constants before, but there are some differences.}

$$\\$$
$\mathbf {Data~ Availability}$
$$\\$$
No data were used to support this study.
$$\\$$
$\mathbf{Conflicts ~of~ Interest}$
$$\\$$
The author(s) declare(s) that there is no conflict of interest regarding the publication of this paper.
$$\\$$
$\mathbf{Funding~ Statement}$
$$\\$$
This work was supported by the National Natural Science Foundation
of P. R. China (Nos. 11971493 and 12071491).


\begin{thebibliography}{99}
 \bibitem{JA3} Alonso, J.;  Martín P.  {\it A counterexample for a conjecture of G. Zbăganu about the Neumann–Jordan constant.} Revue Roumaine Math. Pures Appl. 2006, 51, 135–141.

  \bibitem{Ahmad} Ahmad, A.; Liu,Q.; Li,Y. {\it Geometric Constants in Banach Spaces Related to the Inscribed Quadrilateral of Unit Balls.} Symmetry. 2021, 13, 1294.

\bibitem{JAC} Clarkson J.A.  {\it Uniformly convex spaces.} Trans. Amer. Math. Soc. 1936, 40, 396–414.

\bibitem{Clarkson1937} Clarkson, J.A. {\it The von Neumann–Jordan constant for the Lebesgue space.} Ann. Math. 1937, 38,  114–115.


 \bibitem{Gao1990} Gao,J.; Lau,K.S. {\it On the geometry of spheres in normed linear spaces.} J. Aust. Math. Soc. Ser. A.
1990, 48, 101–112.


\bibitem{Gao2000} Gao, J.  {\it Normal hexagon and more general Banach spaces with uniform normal structure.} Journal of Mathematics.
2000, 20(3), 241-248.


\bibitem{James1964}
James, R. C. {\it Uniformly non-square Banach spaces.} Ann. Math.
	1964, 80, 542–550.

\bibitem{Kato2001}
Kato, M.; Maligranda, L.; Takahashi, Y. {\it On James Jordan–von Neumann constants and the normal
structure coefficient of Banach spaces.} Stud. Math.  2001, 144, 275–295.

\bibitem{LQ} Liu, Q., Zhou, C.; Sarfraz, M.; Li,Y. {\it On New Moduli Related to the Generalization of the Parallelogram Law.} Bull. Malays. Math. Sci. Soc. 2021. https://doi.org/10.1007/s40840-021-01196-7


\bibitem{Mizuguchi2021} Mizuguchi, H. {\it The von Neumann–Jordan and another constants in Radon planes.} Monatsh. Math.
 2021, 195, 307–322.



 \bibitem{ZG}Zbăganu, G. {\it An inequality of M. Rădulescu and S. Rădulescu which characterizes
inner product spaces.} Rev. Roumaine Math. Pures Appl. 2001, 47, 253–257.

\bibitem{Zuo2020} Zuo, Z. F. {\it On James type constants and the normal structure in Banach spaces. Math. Inequal. Appl.} 2020, 23, 341–350.




\end{thebibliography}
\end{document}